\documentstyle[12pt]{article}

  \newcommand{\const}{\rm const}

  \newcommand{\mes}{\rm mes}
  \newcommand{\vraisup}{\rm vraisup}
  \newcommand{\Dom}{\rm  Dom}

\begin{document}

 \begin{center}

 \ {\bf  Sharp Grand Lebesgue Spaces norm estimation }\par

\vspace{4mm}

{\bf for infimal convolution.} \par

\vspace{5mm}

{\bf  M.R.Formica, E.Ostrovsky, and L.Sirota. } \par

 \vspace{5mm}

 \end{center}

 \ Universit\`{a} degli Studi di Napoli Parthenope, via Generale Parisi 13, Palazzo Pacanowsky, 80132,
Napoli, Italy. \\

e-mail: mara.formica@uniparthenope.it \\

\vspace{3mm}

Department of Mathematics and Statistics, Bar-Ilan University,
59200, Ramat Gan, Israel. \\

\vspace{4mm}

e-mail: eugostrovsky@list.ru\\
Department of Mathematics and Statistics, Bar-Ilan University,\\
59200, Ramat Gan, Israel.

\vspace{4mm}

e-mail: sirota3@bezeqint.net \\

\vspace{5mm}

\begin{center}

 \ {\it Abstract.}

 \end{center}

 \ We  derive the non - improvable Grand Lebesgue Space norm estimations for multivariate and multidimensional
operator of infimal convolution. \par

 \vspace{5mm}

  \par

\vspace{5mm}

\section{Notations. Statement of problem. Previous results.}

\vspace{5mm}

 \ Let $ \ f_1(x), f_2(x), \ x \in R^d \ $ be two numerical valued functions. The following function
 $ \ g = f_1 \Box f_2: \ R^d \to R \ $

\begin{equation} \label{inf conv}
g(x) = [f_1 \Box f_2 ](x) \stackrel{def}{=} \inf_{y \in R^d} \left\{ \ f_1(y) + f_2(x - y) \ \right\}
\end{equation}
is named {\it infimal convolution} of $ \ f_1, f_2 \ $ ones. More generally, $ \  g_m[f_1,f_2,\ldots,f_m](x) =  \ $

\begin{equation} \label{multi infimal}
g_m(x) :=   \Box_{j=1}^m f_j (x) := (((f_1 \Box f_2) \Box f_3) \ldots \Box f_m)(x), \ x \in R^d
\end{equation}
or equally

\begin{equation} \label{multi infimal direct}
g_m(x):= \Box_{j=1}^m f_j (x) := \inf   \left\{ \ \sum_{j=1}^m f_j(y_j), \ \sum_{j=1}^m y_j = x, \  \right\}, \ x, y_j \in R^d.
\end{equation}

 \ This operation has many applications in convex analysis, theory of optimization etc., see \cite{Rockafellar},
\cite{Stromberg} and so one. \par

\vspace{4mm}

\ {\bf Our aim in this short report is to estimate the Lebesgue - Riesz norm for this infimal convolution through ones for its components. } \par

\vspace{4mm}

 \ We improve the previous results obtained in particular in an article \cite{Ostrovsky 1}.\par

\vspace{4mm}

 \ Recall that the mentioned Lebesgue - Riesz norm for the (measurable) function $ \ h: R^d \to R \ $ is defined as follows

$$
||h||_p \stackrel{def}{=} \left[ \  \int_{R^d} |h(x)|^p \ dx \ \right]^{1/p},  \ p \ge 1,
$$

$$
||h||_{\infty} \stackrel{def}{=} \vraisup_{x \in R^d} |f(x)|.
$$

 \ As ordinary, $ \ L_p(R^d) = L_p = \{  \  f: R^d \to R, \ ||f||_p < \infty   \ \}. \ $ \par

\vspace{4mm}

\begin{center}

{\sc Grand Lebesgue Spaces (GLS).} \\

\end{center}

 \vspace{4mm}

 \ Let $ \ (a,b) = \const, \ a \ge 1, b \in (a,\infty];  \ $  the case   $ \ b = + \infty \ $  is also not excluded.
\  Let also $ \ p \in[a, b) \ $ or $ \ p \in [a,b]; \ $ evidently, the last case take place iff the value $ \ b \ $ is finite
and greatest than $ \ a. \ $ Let
$ \ \psi_{(a,b)}(p)= \psi=\psi(p)$ be a  function  defined in the domain $ \ (a, b), \ $ not necessarily be finite  in each point
inside of the interval $ (a,b), \ $ such that $ \ \inf \psi(p) > 0. \ $ \par

 \ We can and will suppose without loss of generality  $ \ a = \inf \{p, \psi(p) < \infty \};  \  b = \sup \{p, \psi(p) < \infty\} \ $, so that
$ \Dom[\psi] = [a, b) $  or $ \Dom [\psi] = [a, b], $ of course iff $ \ b < \infty. \ $ \par
 \ When $ \ a > 1, \ b < \infty, \ $ we define formally $ \ \psi(p) = + \infty \ $ for the values $ \ p \notin (a,b). \ $ \par

 \ Denote also

\begin{equation} \label{Psi set}
U\Psi \stackrel{def}{=} \cup_{ a \ge 1, \ b \in (a,\infty)} \Psi_{(a,b)} \ .
\end{equation}

\vspace{4mm}

 \ {\bf Definition 1.1.}  \ The Grand Lebesgue Space (GLS) $ G\psi = G\psi_{(a,b)} $ consists of all
the numerical valued (complex, in general case) measurable functions) $ \{ \  h \ \}, \ h: R^d \to R \ $
 having a finite norm

\begin{equation}\label{GLSnorm}
 ||  h||_{G\psi} \stackrel{def}{=} \sup_{p \in \Dom[\psi] }
\left\{ \frac{\ || \ h \ ||_p}{\psi(p)} \right\}.
\end{equation}

\vspace{3mm}

 \ By definition, $ \ C /\infty := 0. \ $

 \ The function $ \ \psi =\psi(p) \ $ is named ordinary {\it generating function} for this Grand Lebesgue Space  $ \ G\psi. \ $ \par

\vspace{3mm}

  \ These GLS spaces are rearrangement-invariant Banach functional spaces in
the classical sense and were investigated in particular in  many
works, see e.g. \cite{Buld Koz AMS},   \cite{Buldygin-Mushtary-Ostrovsky-Pushalsky},  \cite{caponeformicagiovanonlanal2013},
 \cite{Ermakov etc. 1986},  \cite{Fiorenza-Formica-Gogatishvili-DEA2018},
\cite{fioforgogakoparakoNAtoappear}, \cite{fioformicarakodie2017}, \cite{formicagiovamjom2015},
\cite{Iwaniec1}, \cite{Iwaniec2}, \cite{Kozachenko},  \cite{Ostrovsky 0} - \cite{Ostrovsky 3},
\cite{Samko-Umarkhadzhiev},  \cite{Samko-Umarkhadzhiev-addendum} etc.\par

 \  They were applied in particular in the theory
of probability, especially in the theory of random processes and fields; in the functional analysis - operators theory,
theory of partial differential equations (PDE) etc. \par

\vspace{3mm}

 \ The belonging of the function $ \ f = f(x), \ x \in R^d \ $ is closely related with its tail behavior, where a tail function $ \ T_f(u), \ u \ge 1 \ $
 for the (measurable) function $ \ f \ $  is defined as ordinary

$$
T_f(u) = \mes \ \{x, \ x \in R^d, \ |f(x)| \ge u \ \},  \ u \ge 1,
$$
  and with finiteness of  its norm in appropriate Orlicz - Luxemburg space. A most popular class of these spaces:

$$
\psi_s(p):= p^{1/s}, \ s = \const > 0; \ p \in [1,\infty).
$$

\ Ii is  known for instance that

\begin{equation} \label{tail m}
f \in G\psi_s \ \Leftrightarrow   \ \exists C = C(s) > 0, \hspace{3mm} T_f(u) \le \exp \left(- C(s) u^s \right), \ u \ge 1.
\end{equation}

 \ The value $ \ s = 2 \ $ correspondent to the famous subgaussian case. \par

\vspace{3mm}

 \ Another example (degenerate generating function). Define  for some constant $ \ r \ge 1 \ $

\begin{equation} \label{degener}
\psi_{(r)} (p) = 1, \ p = r; \ \psi_{(r)}(p) = + \infty
\end{equation}
otherwise. The $ \ G\psi_{(r)} \ $ norm of the function $ \ f: R \to R \ $  coincides  with the classical Lebesgue - Riesz one:

$$
||f||\psi_{(r)} = ||f||_r.
$$

\vspace{3mm}

 \ Recall also that the so - called {\it fundamental function} $ \ \phi[G\psi](\delta)\, \ \delta \ge 0 \ $ for these spaces
 is defined as follows

\begin{equation} \label{fund fun}
\phi[G\psi](\delta) \stackrel{def}{=} \sup_{p \in \Dom(f)} \left\{ \ \frac{\delta^{1/p}}{\psi(p)} \ \right\}.
\end{equation}

 \ It is investigated in particular in  \cite{Ostrovsky10}.  \par

\vspace{5mm}

\section{Main result: Lebesgue - Riesz spaces.}

\vspace{5mm}

 \ Let $ \ p \ $ be certain fixed number from the interval $ \ [1,\infty). \ $ Introduce the following
 important variable

\begin{equation}  \label{defin K}
K(d,m,p) \stackrel{def}{=} \sup_{\sum_{j=1}^m ||f_j||_p \in(0, \infty)} \left\{ \ \frac{|| \Box_{j=1}^m f_j||_p \ }{\sum_{j=1}^m||f_j||_p} \ \right\}.
\end{equation}

 \ We set ourselves a goal to calculate the {\it exact} value of these  important for us variable. \par

\vspace{4mm}

 \ {\bf Theorem 2.1.}

\begin{equation} \label{exact Lebesgue}
K(d,m,p) = m^{d/p}, \ m = 1,2,\ldots.
\end{equation}

\vspace{4mm}

 \ {\bf Proof.}  \ {\sc Upper estimate.} \par

 \vspace{3mm}

  \ Introduce the well known {\it dilation} operator

$$
T_{\lambda}[f](x) \stackrel{def}{=} f(\lambda \ x), \ \lambda \in ( 0,\infty), \ f \in L_p(R^d),   \ x \in R^d.
$$

 \ One has

\begin{equation} \label{dilat}
||T_{\lambda} [ f]||_p = \lambda^{-d/p}||f||_p,  \ f \in L_p(R^d).
\end{equation}

  \ Further, one can assume without loss of generality that $ \ f_j(x) \ge 0. \ $
 Denote as above

$$
g_m(x) :=  \left[ \ \Box_{j=1}^m f_j \ \right](x).
$$

 \ Evidently,

$$
g_m(x) \le \sum_{j=1}^m f_j \left(   \ \frac{x}{m}   \  \right) = \sum_{j=1}^m T_{1/m}[f_j](x).
$$
 \ We deduce by virtue of triangle inequality

$$
||g_m||_p \le \sum_{j=1}^m  ||T_{1/m}[f_j]||_p.
$$
 \ It remains to use the relation (\ref{dilat}):

\begin{equation} \label{mdp}
||g_m||_p \le m^{d/p} \ \sum_{j=1}^m ||f_j||_p.
\end{equation}

\vspace{4mm}

 \ {\bf For example:}  a Hilbert norm estimate, i.e. when $ \ p = 2: \ $

$$
||g_m||_2 \le  m^{d/2} \ \sum_{j=1}^m ||f_j||_2.
$$

\vspace{4mm}

 {\sc Lover estimate.} \par

\vspace{4mm}

 \ It is easily to verify that the equality in (\ref{mdp}) is attained if for example
$ \ f_j(x) =  G(x), \  $ where $ \ G(x), \ x \in R^d  \ $ is famous Gaussian density function

$$
G(x) = \exp \left(-||x||^2 \ \right), \ ||x||^2 = (x,x) = \sum_{k=1}^d x_k^2.
$$

  \ In detail, we find solving the following extremal problem

$$
\sum_{k=1}^m G(y_k) \to \min \  / \sum_{k=1}^m y_k = x
$$
that $ \ y_k = x/m. \ $ Therefore

$$
 \ \Box_{k=1}^m G(x) = m G(x/m), \hspace{3mm} || \ \Box_{k=1}^m G(\cdot) \ ||_p  = K(d,m,p) \sum_{k=1}^m ||G(\cdot)||_p.
$$

\vspace{3mm}

 \ This completes the proof of theorem 2.1.\par

 \vspace{5mm}

\section{Main result: Grand Lebesgue Space approach.}

\vspace{5mm}

 \ Let now the function $ \ \psi = \psi_{a,b}(p) \ $  be certain function from the set $  \ \Psi_{(a,b)}, \ 1 \le a < b \le \infty, \ $
see(\ref{Psi set}). We suppose that it may be represented as follows

\begin{equation} \label{psi frac}
\psi(p) = \frac{\nu(p)}{\zeta(p)}, \ p \in (a,b)
\end{equation}
for appropriate such a functions $ \ \nu(\cdot), \ \zeta(\cdot) \ $  belonging at the same set $ \ \Psi_{(a,b)}. \ $ For instance,
$ \ \nu(p) = \psi(p), \ \zeta(p) = 1. \ $ \par

 \ Let once again the function $ \ g_m(x), \ x \in R^d \ $ be defined in (\ref{multi infimal}) or equally in (\ref{multi infimal direct}).\par

 \vspace{4mm}

 \ {\bf Theorem 4.1.} Suppose $ \ \forall j = 1,2,\ldots m \ f_j \in G\psi.  \ $ Then the function $ \ g_m \ $ belongs to the Grand Lebesgue Space
$ \ G\nu \ $ and herewith

\begin{equation} \label{Gpsis estim}
||g_m||G\nu \le \phi_{G\zeta}(m^d) \cdot \sum_{j=1}^m ||f_j||G\psi,
\end{equation}
where (we recall) $ \ \phi_{G\zeta}(\delta) \ $  is the fundamental function of the  Grand Lebesgue space $ \ G \zeta, \ $
see (\ref{fund fun}). \par
 \ Furthermore, the last estimate (\ref{Gpsis estim}) is essentially unimprovable, see an example further. \par

\vspace{4mm}

 \ {\bf Proof.} Let $ \ p \ $ be  an arbitrary number from the segment $ \ (a,b). \ $
 We get from the direct  definition of the norm in the  GLS:

$$
||f_j||_p \le ||f_j||G\psi \cdot \psi(p) =  ||f_j||G\psi  \cdot \frac{\nu(p)}{\zeta(p)}.
$$
 \ One can apply the estimation (\ref{mdp}):

\begin{equation} \label{mdp psi}
||g_m||_p \le m^{d/p} \cdot \frac{\nu(p)}{\zeta(p)} \cdot \sum_{j=1}^m ||f_j||G\psi,
\end{equation}
following

\begin{equation} \label{mdp nu zeta}
\frac{||g_m||_p}{ \nu(p) } \le \frac{m^{d/p}}{ \zeta(p) } \cdot \sum_{j=1}^m ||f_j||G\psi, \ p \in (a,b).
\end{equation}

 \ It remains to take the $ \ \sup \ $ over $ \ p \in (a,b): \ $

\begin{equation} \label{sup p}
||g_m||G\nu \le \phi_{G\zeta}(m^d) \cdot \sum_{j=1}^m ||f_j||G\psi,
\end{equation}
Q.E.D.\par

\vspace{4mm}

 \ {\bf An example:}  choose $ \  a = 1, b \in (a,\infty), \ \nu(p) = \psi(p),  \  $
hence $ \ \zeta(p) = 1 \ $ and following $ \  \phi_{G\zeta} (m^d) = m^d. \ $ We deduce

$$
||g_m||G\psi \le m^d  \cdot \sum_{j=1}^m ||f_j||G\psi.
$$

 \vspace{3mm}

  \ Note in conclusion that the last estimate is essentially non - improvable, for
instance,   when $ \ \psi(p) = \psi_{(r)}(p), \ r \in (1,b). \ $ \par

\vspace{6mm}

\vspace{0.5cm} \emph{Acknowledgement.} {\footnotesize The first
author has been partially supported by the Gruppo Nazionale per
l'Analisi Matematica, la Probabilit\`a e le loro Applicazioni
(GNAMPA) of the Istituto Nazionale di Alta Matematica (INdAM) and by
Universit\`a degli Studi di Napoli Parthenope through the project
\lq\lq sostegno alla Ricerca individuale\rq\rq (triennio 2015 - 2017)}.\par

\end{document}